\documentclass[12pt]{article}

\usepackage{amsmath,amsthm}
\usepackage{amssymb}
\usepackage{color}
\usepackage{xspace}
\usepackage[colorlinks=true,
linkcolor=green,
filecolor=brown,
citecolor=green]{hyperref}
\def\red{\textcolor{red} }
\def\blue{\textcolor{blue} }

\def\a{\ensuremath{\mathcal A}\xspace}
\def\b{\ensuremath{\mathcal B}\xspace}
\def\ib{\ensuremath{\mathcal {IB}}\xspace}
\def\db{\ensuremath{\mathcal {DB}}\xspace}
\def\c{\ensuremath{\mathcal C}\xspace}
\def\ip{\ensuremath{\mathcal IP}\xspace}

\begin{document}
\newtheorem{theorem}{Theorem}
\newtheorem{defn}[theorem]{Definition}
\newtheorem{lemma}[theorem]{Lemma}
\newtheorem{prop}[theorem]{Proposition}
\newtheorem{cor}[theorem]{Corollary}
\begin{center}
{\Large
On Permutations Avoiding the Dashed Patterns 32-41 and 41-32                           \\ 
}

\vspace{10mm}
David Callan  \\
\noindent {\small Department of Statistics, 
University of Wisconsin-Madison,  Madison, WI \ 53706}  \\
{\bf callan@stat.wisc.edu}  \\
May 7, 2014
\end{center}

\vspace{2mm}

\begin{abstract}
We show that  permutations of size $n$ avoiding both of the dashed patterns 32-41 and 41-32 
are equinumerous with indecomposable set partitions of size $n+1$, and deduce a related result.
\end{abstract}

\vspace{1mm}

\section{Introduction} 

\vspace*{-4mm}
 
Permutations avoiding one or two 4-letter (classical) patterns have been studied in detail 
\cite{wikiPattern}, as have permutations avoiding one or two 3-letter dashed 
patterns \cite{pairDashed}. Baxter and Pudwell \cite{vincularSchemes} apply 
enumeration schemes to count 
avoiders for one 4-letter dashed pattern but only isolated results are known for more than one 
4-letter dashed pattern, \cite{sergi07,claesson2009} for example. Here we consider the 
pair of patterns (32-41,\:41-32), whose avoiders are amenable to exact counting. 

Recall that a nonempty set partition of $[n]$ can be maximally decomposed into partitions of 
disjoint subintervals covering $[n]$, called its components, and a partition with 
just one component is \emph{indecomposable}. 
Thus the partition 2/13/4/7/568 has 3 components: 2/13,\ 4,\ 7/568, and, after standardizing 
(replace smallest entry by 1, next smallest by 2, and so on), the 
partitions 2/13,\ 1,\ 3/124 are indecomposable. \textbf{We will always write set partitions, 
as here, in the canonical form of increasing blocks and blocks arranged in increasing order of last entry.}
 
Our main result is a bijection from (32-41,\:41-32)-avoiding permutations of $[n]$ to 
indecomposable partitions of $[n+1]$, counted by sequence 
\htmladdnormallink{A074664}{http://oeis.org/A074664} in 
\htmladdnormallink{OEIS}{http://oeis.org/Seis.html}.
The bijection sends the number of increasing runs in the permutation to
the number of blocks in the partition. Indeed, it sends the largest entries of these 
runs save for the last run to the largest entries of the blocks save for the block containing $n+1$.

We use the term \emph{descent} (in a permutation) to mean a pair of consecutive 
entries $ba$ with $b>a$. 
We call $b$ the \emph{descent initiator} and $a$ the \emph{descent terminator}.

In the patterns 32-41 and 41-32 under consideration, the 32 and 41 both represent descents in the permutation. Defining two (distinct) descents $ba$ and $dc$ in a permutation to be \emph{nested} if one of the intervals $[a,b],\:[c,d]$ is contained in the other, we see that (32-41,\:41-32)-avoiders can be alternatively described as permutations that avoid nested descents. We also have the following useful characterization.
\begin{prop}\label{characterization}
A permutation avoids the patterns $32$-$41$ and $41$-$32$ iff its descent initiators are order isomorphic to its its descent terminators $($both taken left to right\,$)$. \qed
\end{prop}
For example, the 4 descents of the avoider 1\,3\,8\,5\,2\,6\,9\,7\,4 begin with 8,5,9,7 and end with 5,2,7,4 respectively and both 
lists are order isomorphic to 3,1,4,2.

\section{The bijection from avoiders to indecomposables}\label{bigbij}

\vspace*{-4mm}

First, the identity permutation $12\cdots n$ corresponds to the one-block partition of $[n+1]$. 
Now let  $\a_{n}^{*}$ denote the set of non-identity (32-41,\:41-32)-avoiders on $[n]$, and let $\ip_{n+1}^{*}$ 
denote the set of indecomposable partitions of $[n+1]$ with more than one block. (The asterisk ${}^{*}$ 
is a reminder of the omitted element.)
We exhibit the bijection from  $\a_{n}^{*}$  to $\ip_{n+1}^{*}$ as a mapping of 5 (easy) steps to 
produce a partition of $[n+1]$---the 5\emph{-step map}---followed by a tweaking, if necessary, 
to make it indecomposable---the \emph{tweaking map}.

So suppose given $p \in \a_n$ and let us
use $p = 2\ 9\ 5\ 1\ 4\ 10\ 12\ 7\ 8\ 3\ 6\ 15\ 17\ 13\ 14\ 16\ 11$ as a working example with $n=17$.

1. Split $p$ into its increasing runs and append $n+1$ to the last run. 
\[
2\ 9\ / \ 5\ / \ 1\ 4\ 10\ 12\ / \ 7\ 8\ / \ 3\ 6\ 15\ 17\ / \ 13\ 14\ 16\ / \ 11\ 18
\]

\vspace*{-3mm}

2. Call the runs that now end with a left-to-right maximum of $p$ the LRMax runs. Thus the 
last run, ending with $n+1$, is not an LRMax run. Highlight the last entry of each LRMax run (with an enclosing box) 
and underline the remaining entries.
\[
\underline{2\vphantom{y}}\ \framebox{9}\ / \ 5\ / \ \underline{1\ 4\ 10\vphantom{y}}\ \framebox{12}\ / \ 7\ 8\ / \ \underline{3\ 6\ 15\vphantom{y}}\ \framebox{17}\ / \ 13\ 14\ 16\ / \ 11\ 18
\]

\vspace*{-2mm}

%
%
%

3. The first run is rather special: it is an LRMax run since by assumption
$p$ has more than one run and so the first run is not the last run, and it is the only LRMax run 
that may lack underlined entries (which occurs when $p$ starts with a descent). 
Color its underlined entries (if any) red. 
For all later LRMax runs,
compare each underlined entry with the preceding boxed entry. If it is smaller than this 
boxed entry, color it blue; if larger, color it red.
\[
\red{\underline{\textbf{2}\vphantom{y}}}\ \framebox{9}\ / \ 5\ / \ \blue{\underline{\textbf{1\ 4}\vphantom{y}}}\ 
\red{\underline{\textbf{10}\vphantom{y}}}\ \framebox{12}\ / \ 7\ 8\ / \ \blue{\underline{\textbf{3\ 6}\vphantom{y}}}\ 
\red{\underline{\textbf{15}\vphantom{y}}}\ \framebox{17}\ / \ 13\ 14\ 16\ / \ 11\ 18
\]
Thus 1, 4 are blue because they are $<9$, and 3, 6 are blue because they are $<12$. 
It is clear that all but the first LRMax run will contain at least one blue entry 
while there may be no red entries.

4. Transfer all red entries to the last LRMax run, and transfer the blue entries in each LRMax run 
to its left neighbor LRMax run, and leave the non-LRMax runs intact.
\[
\blue{\underline{\textbf{1\ 4}\vphantom{y}}}
\ \framebox{9}\ / \ 5\ / \ \ \blue{\underline{\textbf{3\ 6}\vphantom{y}}}\ \framebox{12}\ / \ 7\ 8\ / \ 
\red{\underline{\textbf{2\ 10\ 15}\vphantom{y}}}\ \framebox{17}\ / \ 13\ 14\ 16\ / \ 11\ 18
\]

\vspace*{-2mm}

5. Each segment is still increasing, and a descent still occurs passing from each segment to the next.
Now forget the colors, underlines, and boxes, 
and arrange the runs in increasing order of last entry.
\[
5\ / \ 7\ 8\ / \ 1\ 4\ 9\ / \ 3\ 6\ 12\ / \ 13\ 14\ 16\ / \ 2\ 10\ 15\ 17\ / \ 11\ 18
\]

\noindent The result of this 5-step process will be  a partition of $[n+1]$ in canonical form. 
If the partition happens to be indecomposable 
(as here), it is the desired indecomposable partition of $[n+1]$. It may, however, be decomposable.
So a further tweaking is required. 

We define the tweaking map to be the identity on indecomposable partitions and 
as follows otherwise. 
With $p=4\ 3\ 1\ 2\ 7\ 6\ 5\ 8\ 10\ 9$, for example, the 5-step map yields 
\[
3\ / \ 1 \ 2\ 4\ / \ 6 \ / \ 5 \ 7\ / \ 8 \ 10 \ / \ 9 \ 11\, ,
\]
a decomposable partition with 3 components: $3\ / \ 1 \ 2\ 4\, ,\ \ 6 \ / \ 5 \ 7\, ,\ \  8 \ 10 \ / \ 9 \ 11 $. 
In this case, simply transfer all entries of the last block of the first component, except for its last entry, to the penultimate (next-to-last) block. Here we get
\[
3\ / \ 4\ / \ 6 \ / \ 5 \ 7\ / \ 1 \ 2 \ 8 \ 10 \ / \ 9 \ 11\, .
\]
The resulting partition will clearly be indecomposable since the last block, containing $n+1$, 
will always contain more than one entry and never consist of consecutive integers.

In the next section, we will see that this mapping (5-step map followed by the tweaking map) is 
a bijection by characterizing the image of the 5-step map and then exhibiting the inverse of both maps.

\section{The inverse mapping}

\vspace*{-4mm}

Here, we show our mapping from  $\a_{n}^{*}$ to $\ip_{n+1}^{*}$ has an inverse. 
To facilitate the description, we make some definitions. 
A \emph{singleton} block in a set partition is a block consisting of only one entry, 
and a \emph{big} block 
is a block consisting of more than one entry. Say that an entry $a$ or a singleton block $\{a\}$ is 
\emph{straddled} by a big block $B$ iff $\min B < a < \max B$.

\begin{lemma}\label{singleton} 
In an indecomposable partition $\pi$ of $[n+1]$ with $n\ge 1$, every singleton block is straddled by a big block.
\end{lemma}

\vspace*{-4mm}

Proof. If the singleton block $\{a\}$ is not straddled by a big block, then $\{a\}$ forms a one-block component of $\pi$. \qed

\begin{lemma}\label{image5} 
A partition $\pi$ of $[n+1]$ obtained from the $5$-step map applied to a non-identity 
permutation $p$ of $[n]$ has the following two properties.

$($i$\,)$ The last block of $\pi$ does not consist of consecutive integers.

$($ii$\,)$ Every singleton block of $\pi$ is straddled by a non-penultimate big block.
\end{lemma}

\vspace*{-3mm}

Proof.  The last block of $\pi$ contains $n+1$ and so was obtained by appending $n+1$ to the 
last run of $p$ and (i) follows. 
As for (ii), suppose $\{a\}$ is a singleton run in $p$. 

If $a$ is the first entry of $p$, then either $a=n$ and $\{a\}$
is a singleton block of $\pi$ that is certainly straddled by the \emph{last} block of $\pi$, or $a<n$ 
in which case $\{a\}$ will not form a singleton block of $\pi$ after the transfer of blue entries.

Otherwise $a$ both starts and ends a descent in $p$ and has an 
LRMax run somewhere to its left in $p$. The run $\{a\}$ will not be disturbed by 
the transfer of colored entries, and will be straddled by the first LRMax 
run to its left after the transfer of blue entries. The only singleton block that can possibly 
be newly-introduced by the transfer of colored entries is one consisting of the largest entry of the 
last LRMax run of step 2, and this largest entry is straddled by the block containing $n+1$. \qed
\begin{defn}\label{defn}
Let $\b_{n+1}$ denote the set of partitions of $[n+1]$ with properties $($i$\,)$ and $($ii\,$)$ of Lemma \ref{image5}.
\end{defn}

\vspace*{-3mm}
 
Split $\b_{n+1}$ into its indecomposable and decomposable elements, denoted $\ib_{n+1}$ and $\db_{n+1}$ respectively. 
Thus $\b_{n+1} = \ib_{n+1} \sqcup \db_{n+1}$ (where $\sqcup$ denotes disjoint union).
Also, split $\ip_{n+1}^{*}$ into two subsets, those in which every singleton block is straddled by a 
non-penultimate big block and the complementary set.
The first subset coincides with $\ib_{n+1}$ because in an indecomposable partition of $[n+1]$ 
with more than one block, the last block, containing $n+1$, cannot consist of 
consecutive integers for otherwise the 
partition would be decomposable. We denote the second subset $\c_{n+1}$.  
Thus $\ip_{n+1}^{*} =  \ib_{n+1} \sqcup \c_{n+1}$.
  
\begin{prop} 
The tweaking map is a bijection from $\b_{n+1}$ onto $\ip_{n+1}^{*}$.
\end{prop}

\vspace*{-4mm}
  
Proof. Here is the inverse of the tweaking map. Let $\pi \in \ip_{n+1}^{*}$. 
If $\pi \notin \c_{n+1}$, fix it to obtain a partition in $\ib_{n+1}$. 
Now suppose $\pi \in \c_{n+1}$. By definition, $\pi$ contains a singleton block $\{a\}$ 
not straddled by any non-penultimate big block. 
Take the largest such $a$. Since $\pi$ is indecomposable, Lemma \ref{singleton} 
implies that its penultimate block is big and straddles $a$, hence contains entries $<a$. 
Transfer all such entries to the block containing $a$ to obtain a partition in $\db_{n+1}$. \qed

\begin{prop} 
The 5-step map is a bijection from the nonidentity $(32$-$41,\:41$-$32)$-avoiding permutations of $[n]$ onto $\b_{n+1}$.
\end{prop}

\vspace*{-4mm}

Proof. 
It is clear from the definition of the 5-step map that the descent initiators in 
the permutation 
can be recovered from the partition 
as the largest entries in the non-last blocks, and the descent 
terminators can be recovered as the smallest entries  
in the non-penultimate blocks. The characterizing order-isomorphic property of the 
descent initiators and terminators (Proposition \ref{characterization})
can then be used to recover the entire avoider as follows, with
\[
5\ / \ 7\ 8\ / \ 1\ 4\ 9\ / \ 3\ 6\ 12\ / \ 13\ 14\ 16\ / \ 2\ 10\ 15\ 17\ / \ 11\ 18
\]
(looks familiar?) as a working example.
First, we retrieve the ordering of the blocks obtained in step 4 (the rest is easy, and left to the reader).
Extract the descent initiators and terminators and list them in increasing order
\[
\begin{array}{lrrrrrr}
\textrm{descent initiators} & 5 & 8 & 9 & 12 & 16 & 17  \\
\textrm{descent terminators} & 1 & 3 & 5 & 7 & 11 & 13
\end{array}
\]
The largest descent initiator, DI for short, is necessarily an LRMax.
So, in step 4, the runs after the run ending at the largest DI 
are unchanged from step 1. 
Hence, the block starting with the largest descent terminator, DT for short, (here, 13) 
immediately follows the block ending at the largest descent initiator (17):
\[
 \ 2\ 10\ 15\ 17\ / \  13\ 14\ 16
\]
The last entry of this block (16) is the 5th largest DI, and so this block is followed by 
the block starting with the the 5th largest DT (11), and we have
\[
\ 2\ 10\ 15\ 17\ / \  13\ 14\ 16\ / \ 11\ 18
\]
Repeat until the block containing $n+1$ is reached. 
(It always will be, and in the example, we're already there.)   
We have now reconstructed a terminal segment.

Start afresh with the largest DI not yet treated (12) and repeat the preceding process until 
we arrive at the DT in the block of this largest DI. We get 
\[
\ 3\ 6\ 12\ / \ 7\ 8\, ,
\]
stopping here because the DT that is order-isomorphic to 8 is 3, which is already in the block of 12. 
Continue in this way to retrieve the full listing in step 4.

There is a more compact description of the entire construction. Let $k$ be the number of blocks. 
Define lists $A=(a_i)_{i=1}^{k}$ and $B=(b_i)_{i=1}^{k}$ by taking $A$ to be the last entries 
of the blocks listed in increasing order and $B$ to be the first entries of the blocks listed 
in increasing order except that the first entry of penultimate block is listed last.
\[
\begin{array}{lrrrrrrr}
\textrm{A\, :\ } & 5 & 8 & 9 & 12 & 16 & 17 & 18  \\
\textrm{B\, :\ } & 1 & 3 & 5 & 7 & 11 & 13 & 2
\end{array}
\]
Define $f: A \to A$ by taking $f(a_{i})$ to be the last entry of the block whose first 
entry is $b_i$. Clearly, $a_k=n+1$, $f(a_k)=a_{k-1}$, and $f$ is a permutation of $A$. Split $f$ into 
disjoint cycles so that each cycle starts with its largest entry except that the cycle containing $a_k=n+1$ starts with $a_{k-1}$, and the cycles are arranged in order of increasing first entries,
\[
(9,\,5) \ (12,\, 8) \ (17,\, 16,\, 18)\, ,
\]
and then erase the parentheses to get a list $C$
\[
C = (9,\,5,\,12,\, 8,\, 17,\, 16,\, 18). 
\]
Now arrange the blocks so their last entries are in the order of $C$ to obtain the 
desired list of segments in step 4.

The construction guarantees that the last entries of the non-last segments are order 
isomorphic to the first entries of the non-first segments. That the segments 
correspond (under step 1 $\leftrightarrow$ step 4) to the maximal 
increasing runs of a permutation is a consequence of the following lemma, taking the $a_i$ 
as the last entries of the blocks and the $b_i$ as the first entries.
\begin{lemma} 
Suppose $(a_i)_{i=1}^{k}$ and $(b_i)_{i=1}^{k}$ are sequences of positive integers with 
the $a_i$ increasing, the $b_i$ distinct, $a_i \ge b_i$ for all $i$, and $a_{k-1}>b_{k}$. 
Suppose also $($the straddling hypothesis\,$)$ that whenever $a_i=b_i$, there exists $j\ne k-1$ 
such that the interval $[b_j,a_j]$ contains $a_i$.

Then, for $1\le i \le k-1$, we have $a_{i}>b_{i}^{*}$ where $(b_{1}^{*}, \dots, b_{k-1}^{*})$ 
is the increasing rearrangement of $(b_1,b_2, \dots, b_{k-2},b_{k})$. \qed
\end{lemma}

\section{A related result}

\vspace*{-4mm}

\begin{prop}
$(32$-$41,\:41$-$32)$-avoiders on $[n]$ that end with $1$ are counted by the Bell number $B_{n-1}$.

\end{prop}
Proof. Here is a bijection from these avoiders to (unrestricted) set partitions of $[n-1]$, 
illustrated using the example 2\:7\:3\:9\:5\:6\:8\:4\:1 with $n=9$. 
Delete the last entry and standardize:
\[
1\:6\:2\:8\:4\:5\:7\:3\, .
\]
Apply the bijection of section \ref{bigbij} from avoiders to indecomposable partitions:
\[
2\:6\ / \ 4\:5\:7\ / \ 1\:8\ / \ 3\:9\, .
\]
Delete the last entry of the last block and then rotate the blocks to the right,
\[
3\ / \ 2\:6\ / \ 4\:5\:7\ / \ 1\:8\, ,
\]
to obtain a partition of $[n-1]$ in canonical form.

We leave the reader to verify that these steps are reversible and that the reverse steps  
produce an avoider from an arbitrary partition of $[n-1]$. \qed


\begin{thebibliography}{99}

\bibitem{wikiPattern} 
\htmladdnormallink{Enumerations of specific permutation classes}{http://en.wikipedia.org/wiki/Enumerations_of_specific_permutation_classes
}, Wikipedia.
 

\bibitem{pairDashed} Claesson, Anders and Mansour, Toufik, Enumerating permutations avoiding a pair of Babson-Steingr\'{i}msson patterns,
\emph{Ars Combin.} \textbf{77} (2005), 17--31.


\bibitem{vincularSchemes}
Baxter, Andrew M. and Pudwell, Lara K., Enumeration schemes for vincular patterns,
\emph{Discrete Math.} \textbf{312} (2012), no. 10, 1699--1712.


\bibitem{sergi07} Sergi Elizalde, Generating Trees for Permutations Avoiding Generalized
Patterns, \emph{Annals of Combinatorics} \textbf{11} (2007) 435--458.

 
\bibitem{claesson2009} Anders Claesson, Sergey Kitaev, Einar Steingr'msson, Decompositions and statistics for $\beta$(1, 0)-trees and nonseparable permutations, \emph{Advances in Applied Mathematics} \textbf(42) (2009) 313Ð-328.
 

\end{thebibliography}
\end{document}